\documentclass[a4paper,12pt]{article}
\usepackage[english]{babel}
\usepackage[T2A]{fontenc}
\usepackage[cp1251]{inputenc}
\usepackage[english]{babel}
\usepackage{amsthm}
\usepackage[tbtags]{amsmath}
\usepackage{amsfonts,amssymb}
\begin{document}

\newtheorem{theorem}{Theorem}
\newtheorem{lemma}{Lemma}
\newtheorem{proposition}{Proposition}
\newtheorem{Cor}{Corollary}

\begin{center}
{\bf Automorphisms of Formal Matrix Rings}
\end{center}
\hfill \textbf{Piotr Krylov}, 

\hfill National Research Tomsk State University, 

\hfill krylov@math.tsu.ru

\hfill \textbf{Askar Tuganbaev}

\hfill National Research University <<MPEI>>, 

\hfill tuganbaev@gmail.com

\textbf{Abstract.} We study automorphism groups of formal matrix algebras. We also consider automorphisms of ordinary matrix algebras (in particular, triangular matrix algebras).


\textbf{MSC2010 database 16S50; 16D10}

\tableofcontents

\section{Introduction}\label{section1} 

Automorphisms and isomorphisms of various matrix rings are studied in many papers; for example, see \cite{AbyT15}, \cite{AnhW11}, \cite{AnhW13}, \cite{HaeH00}, \cite{Isa80}, \cite{Jon91}, \cite{Jon95}, \cite{Kez90}, \cite{KryN18}, \cite{KryN18b}, \cite{KryT21}, \cite{Lev75}, \cite{LiW12}, \cite{Tap15}, \cite{Tap17}, \cite{Tap18}, \cite{XiaW10}, \cite{XiaW14}.

Some other mappings of matrix rings were also studied; in particular, commuting and centralizing mappings were studied (for example, see \cite{LiW12}, \cite{LiWF18}, \cite{XiaW10}, \cite{XiaW14}).

The author's work \cite{KryT21} is devoted to automorphisms and homomorphsms of formal matrix algebras. First, the authors consider automorphisms of the algebra $S=L\oplus M$, where $L$ is some subalgebra and $M$ is a nilpotent ideal. In such a case, one says that $S$ is a \textsf{splitting extension} of the ideal $M$ by the subalgebra $L$. In \cite{KryT21}, results on the group $\text{Aut}\,S$ are applied to formal matrix algebras. At the same time, some assertions are not given in full generality. 

In the given paper, we significantly amplify some results of \cite{KryT21} and give many new results on automorphisms of formal matrix rings. We also study automorphisms of ordinary matrix rings (in particular, triangular matrix rings).

We note that papers \cite{AnhW11} and \cite{AnhW13} contain interesting results and methods of searching for automorphisms  of formal triangular matrix rings. We took and used some ideas from these papers \cite{AnhW11} and \cite{AnhW13}.

We consider only associative rings which are unital algebras over some commutative unital ring $T$. However, the ring $T$ itself is clearly almost non-existent. We sometimes write <<an algebra>>, we sometimes write <<a ring>>.

Let $K$ be some algebra. Then $\text{Aut}\,K$ is the automorphism group of $K$, $\text{In(Aut}\,K)$ is the subgroup of inner automorphisms of $K$, and $\text{Out}\,K$ is the group of outer automorphisms of $K$, i.e., the factor group $\text{Aut}\,K/\text{In(Aut}\,K)$.

If $S$ is a ring, then $U(S)$ is the group of invertible elements and $P(S)$ the prime radical of the ring $S$. For an $S$-$S$-bimodule $M$, we denote by $\text{Aut}_SM$ the automorphism group of $M$.

Let $R$, $S$ be two rings, $A$ be an $R$-$S$-bimodule and let $\alpha$, $\gamma$ be automorphisms of rings $R$ and $S$, respectively.

We can define a new bimodule structure on $A$ by setting
$$
x\circ a=\alpha(x)a,\; a\circ y=a\gamma(y) \text{ for all } x\in R,\, y\in S,\, a\in A.
$$
Usually, this bimodule is denoted by $_{\alpha}A_{\gamma}$ and the initial bimodule can be denoted by $_1A_1$.

The semidirect product of groups $A$ and $B$ is denoted by $A\leftthreetimes B$. This designation has a catchy character but it is convenient. The relation $G\cong A\leftthreetimes B$ implies that the group $G$ contains a normal subgroup $H$ and the subgroup $E$ such that
$$
G=H\cdot E,\; H\cap E=\langle e\rangle,\; A\cong H,\;
B\cong E\cong G/H.
$$.

\section{Group $\text{Aut}\,K$ for Formal Matrix Rings $K$ with zero trace ideals}\label{section2}

The book of authors \cite{KryT17} is devoted to formal matrix rings and formal matrix algebras. We can say on formal matrix algebras, as well.

For ease of reading, we briefly recall some of the material from \cite{KryT17} and \cite{KryT21}.

We fix a positive integer $n\ge 2$. Let $R_1,\ldots,R_n$ be rings and let $M_{ij}$ be $R_i$-$R_j$-bimodules with $M_{ii}=R_i$, $i,j=1,\ldots,n$. Let's assume that for any subscripts $i,j,k=1,\ldots,n$ such that $i\ne j$ and $j\ne k$, an $R_i$-$R_k$-bimodule homomorphsm $\varphi_{ijk}\colon M_{ij}\otimes_{R_j}M_{jk}\to M_{ik}$ is defined. We denote by $\varphi_{iik}$ and $\varphi_{ikk}$ canonical isomorphisms
$$
R_{i}\otimes_{R_i}M_{ik}\to M_{ik},\quad M_{ik}\otimes_{R_k}R_k\to M_{ik}
$$
respectively, $i,k=1,\ldots,n$. Instead of $\varphi_{ijk}(a\otimes b)$, we write $ab$. Using these designations, we also assume that $(ab)c=a(bc)$ for all elements $a\in M_{ij}$, $b\in M_{jk}$, $c\in M_{k\ell}$ and subscripts $i,j,k,\ell$.

We denote by $K$ the set of all square matrices $(a_{ij})$ of order $n$ with values in bimodules $M_{ij}$. With respect to standard operations of matrix addition and matrix multiplication, $K$ forms a ring. We can write it in the following form:
$$
K=\begin{pmatrix}
R_1&M_{12}&\ldots&M_{1n}\\
M_{21}&R_{2}&\ldots&M_{2n}\\
\ldots&\ldots&\ldots&\ldots\\
M_{n1}&M_{n2}&\ldots&R_{n}
\end{pmatrix}.
$$
The ring $K$ is called a \textsf{formal} (or \textsf{generalized}) \textsf{matrix ring} of order $n$. If $M_{ij}=0$ for all $i,j$ with $i>j$, then $K$ is a \textsf{formal (upper) triangular matrix ring}.

For every $k=1,\ldots,n$, we set
$$
I_k=\sum\limits_{i\ne k}\text{Im}(\varphi_{kik}),\;\text{or, in other words,}\; I_k=\sum\limits_{i\ne k}M_{ki}M_{ik}.
$$ 
Here $M_{ki}M_{ik}$ is the set of all finite sums of elements of the form $ab$, where $a\in M_{ki}$ and $b\in M_{ik}$. Then $I_k$ is an ideal of the ring $R_k$. One says that $I_1,\ldots,I_n$ are \textsf{trace ideals} of the ring $K$.

As usual, we identify some matrices with the corresponding elements. For example, we can identify the matrix of the form 
$\begin{pmatrix}
r&0\\
0&0
\end{pmatrix}$
with the element $r$ and so on. Similar agreements also apply to  matrix sets.

Let $K$ be some formal matrix algebra. We denote by $L$ the subring of all diagonal matrices and by $M$ the subgroup of all matrices with zeros on the main diagonal. We can write the direct sum $K=L\oplus M$ of Abelian groups. The subgroup $M$ is an ideal if and only if all trace ideals of the ring $K$ are equal to zero. In this case one says that $K$ is a ring \textsf{with zero trace ideals}. These rings include all triangular matrix rings.

Let $K$ be a formal matrix ring with zero trace ideals. We have a splitting extension $K=L\oplus M$, where $M$ is a nilpotent ideal of nilpotence degree $\le n$ and an $L$-$L$-bimodule. In \cite{KryT21}, automorphisms of such rings $K$ are represented by certain matrices of order $2$. It is done as follows. To an arbitrary automorphism $\varphi$ of the algebra $K$, the matrix 
$\begin{pmatrix}
\alpha&\gamma\\
\delta&\beta
\end{pmatrix}$ can be compared in a standard way. Here
$$
\alpha\colon L\to L,\; \beta\colon M\to M,\;
\gamma\colon M\to L\; \delta\colon L\to M
$$
are $T$-module homomorphsms and
$$
\varphi(x+y)=\begin{pmatrix}
\alpha&\gamma\\
\delta&\beta
\end{pmatrix}
\begin{pmatrix}
x\\
y
\end{pmatrix}=
(\alpha(x)+\gamma(y))+(\delta(x)+\beta(y))
$$
for all $x\in L$ and $y\in M$. Similar to \cite{KryT21}, we mainly consider (except for Section 10) only <<triangular>> case; we mean that $\gamma=0$ for any automorphism $\varphi$. In what follows, we do not distinguish the automorphism $\varphi$ and the matrix corresponding to it. For brevity, we sometimes write <<triangular automorphism $\varphi$>> if $\varphi=\begin{pmatrix}
\alpha&0\\
\delta&\beta
\end{pmatrix}$, and a <<diagonal automorphism $\varphi$>> if $\varphi=\begin{pmatrix}
\alpha&0\\
0&\beta
\end{pmatrix}$.

Let $\varphi=\begin{pmatrix}
\alpha&0\\
\delta&\beta
\end{pmatrix}$ be some automorphism of the algebra $K$. In such case $\alpha$, $\delta$ and $\beta$ satisfy the relations of \cite[Section 3]{KryT21}. In particular, $\alpha$ is an automorphism of the algebra $L$ and $\beta$ is an automorphism of the algebra $M$ (as a non-unital algebra). If $M^2=0$, then $\delta$ is a derivation of the algebra $L$ with values in the bimodule $_{\alpha}M_{\alpha}$ and $\beta$ is an isomorphism of $L$-$L$-bimodules $M\to {}_{\alpha}M_{\alpha}$.

We denote by $\text{In}_1(\text{Aut}\,K)$ (resp., $\text{In}_0(\text{Aut}\,K)$) the subgroup of inner automorphisms of the algebra $K$ defined by invertible elements of the form $1+y$, $y\in M$, (resp., invertible elements of the algebra $L$). The first subgroup is normal in $\text{Aut}\,K$ and we have the semidirect decomposition
$$
\text{In(Aut}\,K)=\text{In}_1(\text{Aut}\,K)\leftthreetimes \text{In}_0(\text{Aut}\,K)
$$
(see \cite[Section 4]{KryT21}).

We define a homomorphism and several groups (see \cite[Section 3]{KryT21}).

Let $f\colon \text{Aut}\,K\to \text{Aut}\,L$ be a homomorphism  such that $f(\varphi)=\alpha$ for every automorphism $\varphi=\begin{pmatrix}
\alpha&0\\
\delta&\beta
\end{pmatrix}$. Next, let $\Lambda$ be the subgroup of diagonal automorphisms and let $\Psi$ be the subgroup consisting of automorphisms of the form $\begin{pmatrix}
1&0\\
0&\beta
\end{pmatrix}$. In addition, we denote by $\Omega$ the image of the homomorphsm $f$. Further, let $\Phi$ be the normal subgroup 
$$
\left\{\varphi\in \text{Aut}\,K \,|\, \varphi=\begin{pmatrix}
\alpha&0\\
\delta&\beta
\end{pmatrix},\;\alpha\in\text{In(Aut}\,L)\right\}
$$
of the group $\text{Aut}\,K$. Information on the defined groups is very important to understand the structure of the group $\text{Aut}\,K$. 

Let $e_1,\ldots,e_n$ be the identity elements of rings $R_1,\ldots,R_n$, respectively. We identify them with the corresponding matrix units.

We formulate two conditions for the algebra $K$ (see \cite[Section 9]{KryT21}).

\textbf{(I)} For any $\varphi\in\text{Aut}\,K$, the relation $\varphi M=M$ holds, i.e., any automorphism is triangular.

\textbf{(II)} For any $\varphi\in\text{Aut}\,K$ and every $i=1,\ldots,n$, we have the inclusion $\varphi(e_i)\in e_k+M$ 
for some $k$.

\textsf{Condition \textbf{(II)} implies condition \textbf{(I)}.}

We formulate in detail and in a more complete form the main results of \cite[Sections 8 and 9]{KryT21} about the group $\text{Aut}\,K$, where $K$ is a formal matrix algebra with zero trace ideals.

First, we write several useful relations and isomorphisms:
$$
\Psi\cap \text{Aut}\,K=\Psi\cap \text{In}_0(\text{Aut}\,K)=\Psi_0;\leqno \textbf{(1)} 
$$
$$
\Lambda/(\text{In}_0(\text{Aut}\,K)\cdot \Psi)\cong
\Omega/\text{In(Aut}\,L); \leqno \textbf{(2)}
$$
$$
\Phi/\text{Ker}\,f\cong \text{In}_0(\text{Aut}\,K)/\Psi_0\cong
\text{In(Aut}\,L);\leqno \textbf{(3)}
$$
$$
\Phi/\text{In(Aut}\,K)\cong \Psi/\Psi_0.\leqno \textbf{(4)}
$$
The group $\Psi_0$ is the subgroup of inner automorphisms of the algebra $K$ defined by the central elements of the algebra $L$ (this group is defined in \cite[Section 4]{KryT21}).

In the following theorem, we collect main information about the group $\text{Aut}\,K$.

\textbf{Theorem 2.1.} Let $K$ be a formal matrix algebra with zero trace ideals such that condition \textbf{(I)} holds. Then we have the following assertions.

\textbf{(a)} We have the relations
$$
\text{Aut}\,K=\text{In}_1(\text{Aut}\,K)\leftthreetimes \Lambda; \leqno \textbf{a1)} 
$$
$$
\text{Ker}\,f=\text{In}_1(\text{Aut}\,K)\leftthreetimes \Psi;\leqno \textbf{a2)}
$$
$$
\Phi=\text{In}(\text{Aut}\,K)\cdot \Psi=
\text{In}_1(\text{Aut}\,K)\leftthreetimes 
\text{In}_0(\text{Aut}\,K)\cdot \Psi.\leqno \textbf{a3)}
$$
\textbf{(b)} We have isomorphisms
$$
\text{Aut}\,K/\text{Ker}\,f\cong \Omega\cong \Lambda/\Psi; \leqno \textbf{b1)} 
$$
$$
\text{Aut}\,K/\Phi\cong \Omega/\text{In}(\text{Aut}\,L).\leqno \textbf{b2)}
$$
\textbf{(c)} The group $\text{Out}\,K$ has a normal 
subgroup which is isomorphic to $\Psi/\Psi_0$ and the factor group with respect to the subgroup is isomorphic to $\Omega/\text{In}(\text{Aut}\,L)$.
 
\textbf{(d)} If the relation $\Omega=\text{In}(\text{Aut}\,L)$  holds, then we have that
$$
\text{Aut}\,K=\Phi=\text{In}_1(\text{Aut}\,K)\leftthreetimes 
(\text{In}_0(\text{Aut}\,K)\cdot \Psi); \leqno \textbf{d1)} 
$$
$$
\text{Out}\,K\cong \Psi/\Psi_0.\leqno \textbf{d2)}
$$
\textbf{(e)} If the relation $\Psi=\Psi_0$ holds, then
$$
\Phi=\text{In}(\text{Aut}\,K), \quad 
\text{Out}\,K\cong \Omega/\text{In}(\text{Aut}\,L).
$$
It can be concluded that if we can find the structure of the groups $\Psi$ and $\Omega$, then the structure of the groups $\text{Aut}\,K$ and $\text{Out}\,K$ is known in some way.

In \cite{KryT21}, the authors calculate the groups $\Psi$ and $\Omega$ for an algebra $K$ over commutative indecomposable ring under some conditions. In Sections $7$ and $8$, these results
receive considerable development and the indicated algebras are defined in Section 5.

A ring $R$ is called \textsf{indecomposable} if $1$ is the unique  non-zero central idempotent of $R$.

\textbf{Corollary 2.2.} Let all factor rings $R_1/P(R_1),\ldots, R_n/P(R_n)$ be indecomposable. Then for the algebra $K$ conditions \textbf{(II)} and \textbf{(I)} hold; consequently, we obtain assertions of Theorem 2.1.

As we agreed, we identify the ring $R_i$ with the ring $e_iKe_i$; we also identify the bimodule $M_{ij}$ with the bimodule $e_iKe_j$.

The example from \cite{Jon95}, also given in \cite{KryT21}, says that automorphisms can <<mix>> the rings $R_i$ and bimodules $M_{ij}$.
In \cite{KryT21}, we highlighted some conditions preventing such mixing. Taking Theorem 2.1(1) into account, we can restrict ourselves to diagonal automorphisms.

\textbf{Theorem 2.3.} Let's assume that all factor rings $R_1/P(R_1),\ldots, R_n/P(R_n)$ are indecomposable. Let $\varphi=\begin{pmatrix}\alpha&0\\ 
0&\beta\end{pmatrix}$ be the diagonal automorphism of the algebra $K$. Then the automorphism $\alpha$ of the algebra $L$ permute the rings $R_1,\ldots, R_n$ and the automorphism $\beta$ of the $L$-$L$-bimodule $M$ permute bimodules $M_{ij}$ in accordance with some permutation $\tau$ of degree $n$. In addition, the restriction of $\beta$ to $M_{ij}$ is a bimodule isomorphism  $M_{ij}\to M_{\tau(i)\tau(j)}$ (with respect to ring isomorphisms $\alpha\big|_{R_i}\colon R_i\to R_{\tau(i)}$ and $\alpha\big|_{R_j}\colon R_j\to R_{\tau(j)}$).

\section{Formal Triangular Matrix Rings}\label{section3}

In \cite{KryT21}, formal triangular matrix rings were not specifically considered. They have a certain specificity which allows you to penetrate more deeply into the structure of both the rings themselves and their automorphism groups; see \cite{BirHKP00}.

As stated in Section 1, we recall that all our rings are $T$-algebras.

For the ring $K$, we give one condition for the rings $R_1,\ldots,R_n$ which is weaker than the condition from Corollary 2.2 and Theorem 2.3; this condition guarantees that condition \textbf{(I)} holds.

\textbf{Definition 3.1 \cite{AnhW13}.}
An idempotent $e$ of a ring $R$ is called \textsf{semicentral} if $(1-e)Re=0$.

A ring $R$ is said to be \textsf{strongly indecomposable} if $1$ is its unique non-zero semicentral idempotent.

For a ring $R$, we consider the following conditions.

\textbf{(1)} $R$ is an indecomposable ring.

\textbf{(2)} $R$ is a strongly indecomposable ring.

\textbf{(3)} The factor ring $R/P(R)$ is indecomposable.

\textbf{(4)} For any idempotent $e$ of the ring $R$, the relation $(1-e)Re=0$ implies the relation $eR(1-e)=0$.

There are the following relations between the above conditions.
$$
(3)\,\Rightarrow\,(2)\,\Rightarrow\,(1),\qquad (2)\,\Rightarrow\,(4).
$$
To conditions \textbf{(I)} and \textbf{(II)} from Section $2$, we add another condition  for the formal matrix ring $K$:

\textbf{(III)} Each of the rings $R_1,\ldots,R_n$ satisfies the above condition \textbf{(4)}.

In Section $2$, it is remarked that condition \textbf{(II)} implies condition \textbf{(I)}. We will show below that condition \textbf{(III)} also implies condition \textbf{(I)} in the <<triangular>> case.

We write again the formal triangular matrix ring $K$ of order $n$ in full form:
$$
K=\begin{pmatrix}
R_1&M_{12}&\ldots&M_{1n}\\
0&R_2&\ldots&M_{2n}\\
\ldots&\ldots&\ldots&\ldots\\
0&0&\ldots&R_n
\end{pmatrix}.
$$
\textbf{Proposition 3.2.} Let $K$ be a formal triangular matrix algebra and let all rings $R_1,\ldots,R_n$ satisfy condition \textbf{(4)}. Then $K$ satisfies condition \textbf{(I)}, i.e., any automorphism of the algebra $K$ is triangular.

$\lhd$ Similar to Section $2$, we denote by $e_1,\ldots,e_n$ the identity elements of rings $R_1,\ldots,R_n$, respectively. We write splitting extension $K=L\oplus M$.

Suppose, on the contrary, that there exists an automorphism $\varphi$ of the algebra $K$ which is not triangular. For every $i=1,\ldots,n$, we have
$$
\varphi(e_i)=g_i+y_i, \; \text{where } g_i\in L, \, y_i\in M.
$$
Here $g_1,\ldots,g_n$ is a complete orthogonal system of idempotents in $L$. We write them with respect to the decomposition $L=R_1\oplus\ldots\oplus R_n$:
$$
g_1=g_1^{(1)}+\ldots+g_n^{(1)},
$$
$$
\ldots\ldots\ldots\ldots\ldots\ldots\ldots\ldots \eqno (1)
$$
$$
g_n=g_1^{(n)}+\ldots+g_n^{(n)}.
$$
All first summands in the relations $(1)$ form a complete orthogonal system of idempotents in ring $R_1$ and so on. 
Since $e_iKe_j=0$ for any $i,j$ with $i>j$, we have that $\varphi(e_i)K\varphi(e_j)=0$ and, consequently, $g_iLg_j=0$ for all same $i$ and $j$.

By considering the relation $(1)$, we obtain the following relations:
$$
\left(g_1^{(n)}+g_1^{(n-1)}+\ldots+g_1^{(2)}\right)R_1g_1^{(1)}=0,
$$
$$
\left(g_2^{(n)}+g_2^{(n-1)}+\ldots+g_2^{(3)}\right)R_2\left(g_2^{(1)}+g_2^{(2)}\right)=0,
$$
$$
\ldots\ldots\ldots\ldots\ldots\ldots\ldots\ldots\ldots\ldots \eqno (2)
$$
$$
g_n^{(n)}R_n\left(g_n^{(1)}+g_n^{(2)}+\ldots+g_n^{(n-1)}\right)=0.
$$
By considering conditions on the ring $R_1,\ldots,R_n$, we can write the following relation:
$$
g_1^{(1)}R_1\left(g_1^{(n)}+g_1^{(n-1)}+\ldots+g_1^{(2)}\right)=0,
$$
$$
\left(g_2^{(1)}+g_2^{(2)}\right)R_2\left(g_2^{(n)}+g_2^{(n-1)}+\ldots+g_2^{(3)}\right)=0,
$$
$$
\ldots\ldots\ldots\ldots\ldots\ldots\ldots\ldots\ldots\ldots \eqno (3)
$$
$$
\left(g_n^{(1)}+g_n^{(2)}+\ldots+g_n^{(n-1)}\right)R_ng_n^{(n)}=0.
$$
It follows from the relations $(3)$ we obtain that
$$
g_k^{(i)}R_kg_k^{(j)}=0 \eqno (4)
$$
for all $i,j,k=1,\ldots,n$ with $i<j$.

Since the automorphism $\varphi$ is not triangular, there exist subscripts $i$ and $j$ ($i<j$) such that $\varphi M_{ij}$ is not contained in $M$. Therefore, it follows from the relation $\varphi M_{ij}=(g_i+y_i)(L\oplus M)(g_j+y_j)$ that $g_iLg_j\ne 0$. This implies that
$$
g_1^{(i)}R_1g_1^{(j)}\oplus\ldots\oplus g_n^{(i)}R_ng_n^{(j)}\ne 0.
$$
Therefore, $g_k^{(i)}R_kg_k^{(j)}\ne 0$ for some $k$. This contradicts to $(4)$. Consequently, the assertion of the proposition is true.~$\rhd$

It is quite common when a formal triangular matrix ring satisfies condition \textbf{(III)}. This is confirmed by the following result.

\textbf{Lemma 3.3.} If a ring $S$ is either semiprime or normal (e.g., commutative), or strongly indecomposable, then $S$ satisfies condition \textbf{(4)}.

$\lhd$ For a normal or strongly indecomposable ring, the assertion is obvious.

We assume that the ring $S$ is semiprime but it does not satisfy condition \textbf{(4)}. Then $S$ contains an idempotent $e$ such that $(1-e)Se=0$ but $eS(1-e)\ne 0$. The ring $S$ can be identified with formal matrix ring of the form 
$\begin{pmatrix}
eSe&eS(1-e)\\
0&(1-e)S(1-e)
\end{pmatrix}$.
In this ring, $eS(1-e)$ is a non-zero nilpotent ideal; this is a contradiction. Consequently, $S$ satisfies condition \textbf{(4)}.~$\rhd$

The following result strengthens \cite[Corollary 9.9(2)]{KryT21}.

\textbf{Corollary 3.4.} Let $K$ be a formal triangular matrix algebra such that each of the rings $R_1,\ldots,R_n$ satisfies condition \textbf{(4)}. Then for $K$, we have assertions of Theorem 2.1.

We give a partial case of Corollary 3.4.

\textbf{Corollary 3.5.} Let's assume that $K$ is a formal triangular matrix algebra such that $R_1,\ldots,R_n$ are commutative rings and $\Omega=\langle 1\rangle$. For example, let 
$R_1=\ldots=R_n=T$, where $T$ is a commutative indecomposable ring and $M_{ij}\ne 0$ for all $i,j$ with $i<j$. Then we have relation $\text{Aut}\,K=\text{In}_1(\text{Aut}\,K)\leftthreetimes \Psi$.

$\lhd$ It follows from the relation $\Omega=\langle 1\rangle$ that $\text{Aut}\,K=\text{Ker}\,f$. Therefore, the relations from the corollary follow from Corollary 3.4 and Theorem 2.1.

We pass to a partial case. 
It follows from of assertions 2.1, 3.2 and 3.3 that $\text{Aut}\,K=\text{In}_1(\text{Aut}\,K)\leftthreetimes \Lambda$. Let 
$\varphi=\begin{pmatrix}
\alpha&0\\0&\beta
\end{pmatrix}\in \Lambda$. The automorphism $\beta$ permute bimodules $M_{ij}$ in accordance with some permutation $\tau$ (Theorem 2.3). Since all bimodules $M_{ij}$ are non-zero, we have that $\tau$ is the identity permutation. Therefore, it follows from the relations $R_1=\ldots=R_n=T$ that $\alpha=1$; then we obtain the relation $\Omega=\langle 1\rangle$.~$\rhd$

We specialize Theorem 2.3 to the case $K$ of formal triangular matrix rings.

\textbf{Corollary 3.6.} Let $K$ be a formal triangular matrix algebra such that condition \textbf{(I)} holds and the factor rings $R_1/P(R_1),\ldots,R_n/P(R_n)$ are indecomposable. For example, let the rings $R_1,\ldots,R_n$ be strongly indecomposable. In addition, let $M_{ij}\ne 0$ for all $i,j$ with $i<j$.
Then any diagonal automorphism $\varphi=\begin{pmatrix}\alpha&0\\ 0&\beta\end{pmatrix}$ of the algebra $K$ leaves each of the rings $R_1,\ldots,R_n$ in place and the restriction $\beta\big|_{M_{ij}}$ is an isomorphism of $R_i$-$R_j$-bimodules $M_{ij}\to {}_{R_i}(M_{ij})_{R_j}$ for all $i,j$ with $i<j$.

$\lhd$ Automorphisms $\alpha$ and $\beta$ permute the rings $R_1,\ldots,R_n$ and bimodules $M_{ij}$ in accordance with some permutation $\tau$. Similar to the proof of Corollary 3.5, we obtain that $\tau$ is the identity permutation.~$\rhd$ 

\section{Subgroup $\Psi$ and Inner Automorphisms}\label{section4}

In the beginning of the section, $K$ denotes some formal matrix algebra with zero trace ideals.

In \cite[Section 3]{KryT21}, we formulated the calculation problem of the subgroup $\Psi$. Theorem 2.1 shows an important role of this subgroup in the description problem for the automorphism group of the algebra $K$. We give several remarks on the subgroup $\Psi$. Everything below is true for any algebra $K$ (i.e., it is not assumed that all its automorphisms are triangular).

We take an arbitrary automorphism $\varphi=\begin{pmatrix}1&0\\ 0&\beta\end{pmatrix}\in \Psi$. It is known that $\beta$ is an automorphism of the algebra $M$ (as a non-unital algebra) and an automorphism of the $L$-$L$-bimodule $M$. Conversely, if the mapping $\beta$ satisfies these properties, then $\begin{pmatrix}1&0\\ 0&\beta\end{pmatrix}$ is an automorphism of the algebra $K$ contained in $\Psi$. We clarify this observation as follows.

The automorphism $\beta$ induces the automorphism $\beta_{ij}$ on every $R_i$-$R_j$-bimodule $M_{ij}$. At the same time, for any pairwise distinct subscripts $i,j,k$ and elements $a\in M_{ij}$, $b\in M_{jk}$ the relation
$$
\beta_{ik}(a\cdot b)=\beta_{ij}(a)\cdot \beta_{jk}(b) \eqno (*)
$$
must be carried out. For any two subscripts $i,j$, let we have an automorphism $\beta_{ij}$ of the bimodule $M_{ij}$ such that the relation $(*)$ is true for all values of symbols included in it. By setting
$$
\beta=\sum_{i,j=1,\,i\ne j}^n\beta_{ij},
$$
we pass to an automorphism $\begin{pmatrix}1&0\\ 0&\beta\end{pmatrix}$ which belongs to the subgroup $\Psi$.

We obtain a group embedding
$$
\Psi\to \text{Aut}_LM= \prod_{i,j=1,\,i\ne j}^n\text{Aut}_LM_{ij}
$$
if we assign  the set of restrictions $\beta\big|_{M_{ij}}$ to the automorphism $\begin{pmatrix}1&0\\ 0&\beta\end{pmatrix}$ from the subgroup $\Psi$. If the algebra $K$ satisfies the property $M^2=0$, then the correspondence $\begin{pmatrix}1&0\\ 0&\beta\end{pmatrix}\to \beta$ defines isomorphism groups $\Psi\cong \text{Aut}_LM$.

There is another situation where it is also possible to specify the structure of the subgroup $\Psi$. To reveal this situation, we impose additional restrictions condition on the ring $K$. Namely, for any subscripts $i,j,k$ with $i<j<k$, we set that $R_i$-$R_k$-bimodule homomorphsms $\varphi_{ijk}\colon M_{ij}\otimes M_{jk}\to M_{ik}$, used in the definition of the multiplication in $K$, are isomorphisms.

In Section 2, we agreed to write $ab$ instead of $\varphi_{ijk}(a\otimes b)$; in addition, the symbol $M_{ij}\cdot M_{jk}$ denotes the set of all finite sums of elements of the form $ab$, where $a\in M_{ij}$, $b\in M_{jk}$. In other words, $M_{ij}\cdot M_{jk}$ is the image of the homomorphsm $\varphi_{ijk}$. It is also clear what we mean by the product of several bimodules $M_{ij}$. Thus, for $i<j<k$ or $i>j>k$, we have the relation
$M_{ij}\cdot M_{jk}=M_{ik}$. Also we have the relations
$$
M_{ik}=M_{i,i+1}\cdot M_{i+1,i+2}\cdot\ldots\cdot M_{k-1,k},
$$
$$
M_{ik}=M_{i,i-1}\cdot M_{i-1,i-2}\cdot\ldots\cdot M_{k+1,k},
$$
respectively, for $i<k$ and $i>k$.

Next, we obtain the following. For every $i=1,\ldots,n-1$, let we have an automorphism $\beta_{i,i+1}$ of the bimodule $M_{i,i+1}$. These automorphisms induce the uniquely defined automorphism $\beta_{ik}$ of the bimodule $M_{ik}$ for all $i,k$ with $i<k$. Similarly, the set of automorphisms $\beta_{i,i-1}$ of bimodules $M_{i,i-1}$ for $i=2,\ldots,n$ induces the automorphism $\beta_{ik}$ of the bimodule $M_{ik}$ for all $i,k$ with $i>k$. At the same time, the relation $(*)$ is true. Thus, automorphisms $\beta_{i,i+1}$ ($i=1,\ldots,n-1$) and $\beta_{i,i-1}$ ($i=2,\ldots,n$) induce the uniquely defined automorphism $\beta$ of the algebra $M$ and the $L$-$L$-bimodule $M$. Consequently, the automorphism 
$\begin{pmatrix}1&0\\ 0&\beta\end{pmatrix}$ are contained in the subgroup $\Psi$. We can write the following result.

\textbf{Corollary 4.1.} In the above situation, we have an isomorphism
$$
\Psi\cong \prod_{i=1}^{n-1}\text{Aut}(M_{i,i+1})\times 
\prod_{i=2}^{n}\text{Aut}(M_{i,i-1}).
$$
If $K$ is a triangular matrix ring, then there is no the second factor in the right part.

In the second half of this section, we touch on the following familiar question: when are all automorphisms of the algebra $K$ inner? It was considered in \cite[Section 10]{KryT21} for triangular matrix algebras. Here much depends of the subgroup $\Psi$.

\textsf{Up to the end of this section, we assume that the algebra $K$ of formal matrices with zero trace ideals satisfies condition \textbf{(I)} given in Section 2.}

The following facts follow from Theorem 2.1(c) and the relations $(1)$ before this theorem.

\textbf{Corollary 4.2.}\\
\textbf{1.} Every automorphism of the algebra $K$ is inner if and only if we have the relation
$$
\Psi=\Psi_0 \text{ and } \Omega=\text{In(Aut}\,L).
$$
\textbf{2.} The inclusion $\Psi\subseteq\text{In(Aut}\,K)$ is equivalent to the relation $\Psi=\Psi_0$.

We give several remarks related to the relation $\Psi=\Psi_0$. It is hardly possible to find criteria for the fulfillment of this relation without additional information about rings $R_i$ and bimodules $M_{ij}$.

What can be said about automorphisms from $\Psi_0$? Let an automorphism $\varphi=\begin{pmatrix}1&0\\ 0&\beta\end{pmatrix}$ belong to $\Psi_0$ and be defined by an invertible central matrix $v=\text{diag}(v_1,\ldots,v_n)\in L$, where $L=\oplus_{i=1}^nR_i$. For any distinct subscripts $i,j$ and any $y\in M_{ij}$, we have the relation
$$
\varphi(y)=\beta(y)=v^{-1}yv=v_i^{-1}yv_j; \eqno (1)
$$
that's all we know about $\varphi$.

\textsf{Now we assume that the rings $R_1,\ldots,R_n$ have pairwise isomorphic centers.} We identify these centers and say the <<common center>>. We denote it by $C$.

Let's assume that the automorphism group of every $R_i$-$R_j$-bimodule $M_{ij}$ consists of multiplications by invertible elements from the center $C$, i.e., $\text{Aut}(M_{ij})=U(C)$. We take $\varphi=\begin{pmatrix}1&0\\ 0&\beta\end{pmatrix}\in \Psi$. Let $c_{ij}$ be an invertible element of the ring $C$ such that the relation $\beta(y)=c_{ij}y$ holds for $y\in M_{ij}$. For the automorphism $\varphi$, the relation $(*)$ takes the form
$$
c_{ik}ab=c_{ij}c_{jk}ab \eqno (**)
$$
$$
\text{or } (c_{ik}-c_{ij}c_{jk})M_{ij}M_{jk}=0. \eqno (***)
$$
Thus, we can assign a system of invertible elements $c_{ij}$, $i,j=1,\ldots,n$, to the automorphism $\varphi$, where we assume that $c_{ii}=1$. For these elements, relations $(**)$ and $(***)$ hold.

The above group embedding $\Psi\to \text{Aut}_LM$ turns into an embedding $\Psi\to \prod_{n^2-n}U(S)$. It is difficult to find the image of it. In Section 7, we do this for the matrix ring over a given ring $R$.

\textbf{Proposition 4.3.} Let $K$ be the algebra from Theorem 2.1. In addition, we assume that all rings $R_i$ have a common center $C$ and $cm=mc$ for all $c\in C$ and $m\in M$. Under such assumptions, the relation $\Psi=\Psi_0$ is true if and only if for any automorphism $\varphi=\begin{pmatrix}1&0 \\ 0&\beta\end{pmatrix}\in \Psi$, there exist such invertible elements $c_{ij}\in C$, $i,j=1,\ldots,n$ that $c_{ii}=1$,

\textbf{a)} $\varphi(y)=c_{ij}y$ for any $i,j$ and $y\in M_{ij}$;

\textbf{b)} $c_{ij}\cdot c_{jk}=c_{ik}$ for all $i,j,k$.

$\lhd$ We assume that the relation $\Psi=\Psi_0$ is true and $\varphi=\begin{pmatrix}1&0 \\ 0&\beta\end{pmatrix}\in \Psi$. Continuing the relation \textbf{(1)}, we obtain $\varphi(y)=v_i^{-1}v_jy$. We set
$$
c_{ij}=v_i^{-1}v_j \text{ and } c_{ii}=1 \text{ for all } i,j.
$$
Elements $c_{ij}$ satisfy \textbf{a)} and \textbf{b)}.

Conversely, let for every automorphism $\varphi\in \Psi$, exist elements $c_{ij}$ with properties mentioned in \textbf{a)} and \textbf{b)}. We choose some invertible element $v_1$ in $C$ and we set $v_2=v_1c_{12}$,$\ldots$,$v_n=v_1c_{1n}$. Then $v_i^{-1}v_j=c_{ij}$ for all $i,j$. In addition, conjugation by the matrix $\text{diag}(v_1,\ldots,v_n)$ coincides with the automorphism $\varphi$. Consequently, $\varphi\in \Psi_0$ and $\Psi=\Psi_0$.~$\rhd$

\textbf{Corollary 4.4.} To conditions of Proposition 4.3, we add another condition:\\
$M_{ij}M_{jk}$ is a faithful $C$-module for all $i,j,k$.\\
Then we can exclude item \textbf{b)} of the proposition.

As for the relation $\Omega=\text{In(Aut}\,L)$, Section 8 contains various information about the group $\Omega$ for the ring formal matrices over a given ring $R$.

\section[Formal Matrix Rings over a Given Ring]{Formal Matrix Rings \\ over a Given Ring}\label{section5}

There is an interesting form of formal matrix rings. They are listed in the title of this section. Such rings are considered in the book \cite{KryT17}.

Namely, let $R$ be some ring. If $K$ is a formal matrix ring such that $R_1=\ldots=R_n=R$ and $M_{ij}=R$ for all distinct subscripts $i$ and $j$, then one says that $K$ is a \textsf{formal matrix ring over the ring} $R$ or \textsf{formal matrix ring with values in the ring} $R$.

We denote by $e_{ij}$ ($i,j=1,\ldots,n$) matrix units of the ring $K$. For all values of subscripts $i,j,k$, we have $e_{ij}e_{jk}=s_{ijk}e_k$ for some central elements $s_{ijk}$ of the ring $R$. Under multiplication of matrices $A=(a_{ij})$ and $B=(b_{ij})$ from $K$, we need to take into account the relation
$$
c_{ij}=\sum_{k=1}^ns_{ikj}a_{ik}b_{kj}, \eqno (*)
$$
where $A\cdot B=C=(c_{ij})$. Elements $s_{ijk}$ satisfy identities
$$
s_{iik}=1=s_{ikk},\; s_{ijk}\cdot s_{ik\ell}=s_{ij\ell}\cdot s_{jk\ell}. \eqno (1)
$$ 
Now let $\{s_{ijk}\,|\,i,j,k=1,\ldots,n\}$ be some set of central of elements of the ring $R$ which satisfy identities $(1)$. If we define multiplication of matrices $A=(a_{ij})$ and $B=(b_{ij})$ by the relation $(*)$, then we obtain a formal matrix ring over the ring $R$. Therefore, two given definitions are equivalent.

Let $K$ be some formal matrix ring over the ring $R$ and let $\Sigma =\{s_{ijk}\,|\,i,j,k=1,\ldots,n\}$ be the corresponding system of central elements. The set $\Sigma$ is called a \textsf{multiplier system} and its elements are called \textsf{multipliers} of the ring $K$. Instead of <<multipliers>> also one says <<multiplicative coefficients>>; for example, see \cite{Tap15}. The ring $K$ can be denoted by $M(n,R,\Sigma)$. If all $s_{ijk}$ are equal to $1$, then we obtain an ordinary matrix ring $M(n,R)$.

Let $\tau$ be a permutation of degree $n$. For any matrix $A=(a_{ij})$ of order $n$, we set $\tau A=(a_{\tau(i)\tau(j)})$, i.e., we take the conjugation of the matrix $A$ by the matrix permutation $\tau$. Next, if $\Sigma=\{s_{ijk}\}$ is some multiplier system, then we set $t_{ijk}=s_{\tau(i)\tau(j)\tau(k)}$. Then $\{t_{ijk}\}$ also is a multiplier system, since it satisfies identities $(1)$. We denote it by $\tau\Sigma$. Consequently, there exists a formal matrix ring $M(n,R,\tau\Sigma)$. The rings $M(n,R,\Sigma)$ and $M(n,R,\tau\Sigma)$ are isomorphic under the correspondence $A\to \tau A$.

\textsf{Up to the end of the section, we assume that $K$ is a formal matrix ring over given ring $R$ which is a $T$-algebra. Also we assume that every multiplier $s_{ijk}$ is equal to $1$ or $0$.}

In \cite{KryT21}, it is shown that, under this assumption, there exists a permutation $\tau$ with property that the ring $\tau K$ can be represented as a ring of formal block matrices with zero trace ideals. We briefly recall this material. 

By considering identities $(1)$, it is easy to verify the following lemma.

\textbf{Lemma 5.1.} Let subscripts $i,j,k$ be pairwise distinct. Then for elements $s_{iji}$, $s_{jkj}$ and $s_{kik}$, we have one of the following possibilities.

\textbf{1)} All three elements are equal to $1$.

\textbf{2)} Some two of these three elements are zeros and the third element is $1$.

\textbf{3)} All three elements are zeros.

On the set of integers $\{1,\ldots,n\}$, we define a binary relation $\sim$ by setting $i\sim j$ $\Leftrightarrow$ $s_{iji}$ is equal to $1$.

\textbf{Lemma 5.2.} The relation $\sim$ is an equivalence relation.

The symmetrical matrix $S=(s_{iji})$ is called the \textsf{multiplier matrix} of the ring $K$.

We construct a permutation $\tau$ as follows. To the upper row, we arrange positive integers from $1$ to $n$ in natural order. The bottom row consists of equivalence classes with respect to relation $\sim$ arranged in an arbitrary order. Inside classes, these integers are also arranged in an arbitrary order. Then the main diagonal of the matrix $\tau S$ contains blocks consisting of 1's. There is an one-to-one correspondence between these blocks and equivalence classes with respect to the relation $\sim$. The order of this block is equal to the number of elements of the corresponding equivalence class. In the matrix $\tau S$, all positions outside the considered blocks are occupied by zeros.

As it was mentioned above, the rings $K$ and $\tau K$ are isomorphic under the correspondence $A\to \tau A$, $A\in K$. To simplify the text, we agree that the multiplier matrix $S$ of the ring $K$ already has the above block form. Let the number of blocks on the main diagonal of matrix $S$ be equal $m$.

On the main diagonal of any matrix $A\in K$, we select blocks $A_1,\ldots,A_m$ of the same order and in the same sequence as on the main diagonal of the matrix $S$. For a fixed $\ell$, the blocks $A_{\ell}$ of all matrices in $K$ form the usual matrix ring $M(k_{\ell},R)$ for some $k_{\ell}$. We denote it by $R_{\ell}$. Blocks $A_1,\ldots,A_m$ define an obvious block decomposition of matrices $A$.

The symbol $L$ denotes the direct sum of rings $R_1\oplus\ldots\oplus R_m$. By $M$, we denote the set of all matrices $A\in K$ such that the corresponding blocks $A_1,\ldots,A_m$ consist of zeros. It is clear that $M$ is an $L$-$L$-bimodule.

The decomposition $L=R_1\oplus\ldots\oplus R_m$ induces the block decomposition of every matrix mentioned above. Namely, we write $1=e_1+\ldots+e_m$, where $e_i$ is the identity element of the ring $R_i$. Now we denote by $M_{ij}$ the subbimodule $e_iMe_j$ in $M$. The action of the ring $L$ on the subbimodule $M_{ij}$ coincides with the action of rings $R_i$ and $R_j$ from the left and the right, respectively.
We have a bimodule direct decomposition $M=\oplus_{i,j=1}^mM_{ij}$, where $i\ne j$. Similar to Section $2$, we have the direct sum $K=L\oplus M$.

The ring $K$ is a formal (block) matrix ring constructed from the rings $R_1,\ldots, R_m$ and bimodules $M_{ij}$ in accordance with  procedure given in Section 2; see \cite[Section 2.3]{KryT17}. Basically, we will consider the ring $K$ as a ring of block matrices.

As a ring of block matrices, the algebra $K$ has zero trace ideals. Consequently, we get into the situation of Section 2.

If the factor ring $R/P(R)$ is indecomposable, then all factor rings $R_i/P(R_i)$ ($i=1,\ldots,m$) are indecomposable, as well. By considering Corollary 2.2, we can write such result.

\textbf{Corollary 5.3.} Let the ring $R/P(R)$ be indecomposable. Then  the algebra $K$ satisfies conditions \textbf{(II)} and \textbf{(I)}; consequently, Theorems 2.1 and 2.3 are true for the group $\text{Aut}\,K$.

\textbf{Remark 5.4.} For a commutative ring $R$, the ring $R/P(R)$ is indecomposable if and only if the ring $R$ is indecomposable. Therefore, if $R$ is an indecomposable commutative ring, then Theorems 2.1 and 2.3 are true for the automorphism group of the $R$-algebra $K$.

\section{Case $M^2=0$}\label{section6}

We preserve all designations and agreements of the previous section. Thus, $K$ is a formal matrix algebra over a given ring $R$ and $K=L\oplus M$, where the symbols $L$ and $M$ have the same meaning. We consider $K$ as a ring of formal block matrices in accordance with Section 5. For every algebra $K$ with $M^2=0$, the authors clarify some facts in \cite{KryT21}; additional information about the group $\text{Aut}\,K$ is also obtained in \cite{KryT21}.

We recall this material. First, we briefly repeat some general considerations from \cite[Section 9]{KryT21}. The following fact is true; see \cite[Lemma 9.6]{KryT21}.

\textbf{Lemma 6.1.} Let we have indecomposable rings $R_1,\ldots,R_n$, $L=R_1\oplus\ldots\oplus R_n$, and let $\alpha\in \text{Aut}\,L$. Then for every subscript $i=1,\ldots,n$, there exists a subscript $j$ such that $\alpha R_i=R_j$.

In \cite{KryT21}, based on Lemma 6.1, a certain permutation group $\Sigma$ of degree $n$ is defined; it acts on the ring $L$. At the same time, a permutation $\sigma\in\Sigma$ is identified with the corresponding automorphism $\alpha_{\sigma}$ of the ring $L$. In this article, it is also introduced a normal subgroup $\Gamma$ of automorphisms of the ring $L$, leaving all $R_i$ in place. Then
$$
\Gamma\cap \Sigma=\langle 1\rangle \, \text{ and }
\text{Aut}\,L=\Gamma\leftthreetimes \Sigma.
$$
We return to formal matrix algebras $K$ over $R$, where $R$ is some ring. Similar to Section 2, we write
$$
K=L\oplus M, \text{ where } L=R_1\oplus\ldots\oplus R_m 
$$
and every $R_i$ is an ordinary matrix ring $M(k_i,R)$ for some $k_i\ge 1$.

We impose the same restrictions on the ring $R$ as in Section 5. Namely, we assume that the ring $R/P(R)$ is indecomposable. Then all rings $R_i/P(R_i)$ are indecomposable, as well. Consequently, the ring $K$ satisfies condition \textbf{(II)}. In addition, all rings $R_i$ are indecomposable, as well. Therefore, we can apply Lemma 6.1 to $L$. We can write $\text{Aut}\,L=\Gamma\leftthreetimes \Sigma$, where $\Gamma$ and $\Sigma$ are such subgroups as indicated after Lemma 6.1.

Let $h\colon \text{Aut}\,L\to \Sigma$ be the canonical homomorphsm and $g=hf\colon \text{Aut}\,L\to \Sigma$. Then we have
$$
\text{Ker}\,g=\left\{\varphi=\begin{pmatrix}\alpha&0\\ \delta&\beta\end{pmatrix} \,|\, \alpha R_i=R_i, \, i=1,\ldots,m\right\}.
$$
With the use of relation $(1)$ in this section $5$, it is easy to verify the following fact.

\textbf{Lemma 6.2.} For a given algebra $K$, the relation $M^2=0$ is true if and only if the following condition holds:\\
for any pairwise distinct subscripts $i,j,k$, the relations $s_{iji}=s_{jkj}=s_{kik}=0$ imply the relation $s_{ijk}=0$.

In \cite[Section 3]{KryT21}, a subgroup $\Delta$ is defined. It consists of automorphisms of the form $\begin{pmatrix}1&0\\ \delta&1\end{pmatrix}$. For our algebra $K$, we have relations
$$
\Delta=\text{In}_1(\text{Aut}\,K) \text{ and } \Delta\cong 1+M.
$$
According to Theorem 2.1 we have the relation $\text{Aut}\,K=\Delta\leftthreetimes \Lambda$. We can also write decomposition $\text{Ker}\,g=\Delta\leftthreetimes C$,
where $C$ denotes 
$$
\left\{\begin{pmatrix}
\alpha&0\\0&\beta
\end{pmatrix}\, \Big|\, \alpha R_i=R_i \text{ for all } i=1,\ldots,m\right\}.
$$
By considering Theorem 2.1, we can formulate the following theorem.

\textbf{Theorem 6.3.} Let $R$ be a ring with indecomposable factor ring $R/P(R)$ and let $K$ be a formal matrix algebra over $R$ with $M^2=0$.

\textbf{1.} We have the relations
$$
\text{Aut}\,K=\Delta\leftthreetimes \Lambda,\;
\text{Aut}\,K=\Delta\leftthreetimes C\leftthreetimes \Sigma,\;
\Phi=\Delta\leftthreetimes (\text{In}_0(\text{Aut}\,K))\cdot\Psi.
$$
\textbf{2.} If all automorphisms of every algebra $R_1,\ldots,R_m$ are inner, then we have relations
$$
\text{Aut}\,K=\Delta\leftthreetimes (\text{In}_0(\text{Aut}\,K)\cdot\Psi)\leftthreetimes \Sigma,
$$
$$
\text{Out}\,K\cong \Psi/\Psi_0\leftthreetimes \Sigma.
$$
\textbf{Remarks on item 2:} When the conditions of this item are met, the relation $\text{Ker}\,g=\Phi$ is true and
$$
\text{Aut}\,K=\text{Ker}\,g\leftthreetimes \Sigma=\Phi\leftthreetimes \Sigma=\Delta\leftthreetimes (\text{In}_0(\text{Aut}\,K)\cdot \Psi)\leftthreetimes\Sigma. 
$$
Included in Theorem 6.3, the structure of subgroups $\Delta$, $\text{In}_0(\text{Aut}\,K)$, $\Psi$ and $\Sigma$  is known (see Sections 7 and 8). Therefore, we know the structure of the whole group $\text{Aut}\,K$ from this item. For example, the condition on automorphisms of the algebras $R_1,\ldots,R_m$ is satisfied for a commutative ring $R$ which is a unique factorization domain or a local ring.

\section[Subgroup $\Psi$ and Inner Automorphisms, II]{Subgroup $\Psi$ and\\ Inner Automorphisms, II}\label{section7}

Section 4 contains some information about the subgroup $\Psi$ for the formal matrix algebra with zero trace ideals. In Section 7, we will calculate this subgroup for a formal matrix ring with values in an arbitrary ring $R$. Thus, we continue the line of Sections 5 and 6. At the same time, we develop the results from \cite[Section 13]{KryT21}. We preserve all designations and terms of Sections 5 and 6. In addition, we assume that $C(U(R))$ is the center of the group $U(R)$.

Let $K$ be a formal matrix algebra with values in a ring $R$. According to Section 4, there exists a group embedding
$$
\Psi\to\text{Aut}_LM=\prod_{i,j=1,i\ne j}^m\text{Aut}_LM_{ij}.
$$
By \cite[Proposition 13.2]{KryT21}, automorphisms of the $R_i$-$R_j$-bimodule $M_{ij}$ coincide with multiplications by invertible central elements of the ring $R$. Therefore, we can  write $\text{Aut}_LM_{ij}=C(U(R))$.

Let $\varphi=\begin{pmatrix}1&0\\ 0&\beta\end{pmatrix}\in \Psi$. We have a system of invertible central elements $c_{ij}$ ($i,j=1,\ldots,m$) of $R$ with $c_{ii}=1$, which satisfy to the relations $(**)$ and $(***)$ from Section 4 for all values of subscripts $i,j,k$. We show that for our algebra $K$,
we can limit ourselves in a certain sense to a smaller number of elements $c_{ij}$. And also we will exactly specify the image of the embedding $\Psi\to \prod_{m^2-m}C(U(R))$. For this purpose, we give some argument.

We fix a subscript $i$, where $1\le i\le m-1$. Let $k_i$ be a subscript such that 
$$
i< k_i\le m,\; M_{i,i+1}\cdots M_{k_i-1,k_i}\ne 0
$$
and $k_i$ is the maximal number with such a property (the meaning of the written product of bimodules is explained in Section 4). Then $M_{ij}\cdot M_{jk}\ne 0$ for any $k$ and $j$ such that $i<k\le k_i$ and $i<j<k$. It follows from the relation $(***)$ from Section 4 that $c_{ik}=c_{ij}c_{jk}$ (note that multipliers $s_{ijk}$ only take values $0$ or $1$).

If there are indices $\ell$ such that $k_i<\ell\le m$, then the elements $c_{i\ell}$ and $c_{ij}c_{j\ell}$ ($i<j< \ell$) may not be related in any way, since $M_{ij}M_{j\ell}=0$. And accordingly, the element $c_{i\ell}$ does not depend on the elements $c_{i,i+1},\ldots, c_{\ell-1,\ell}$.

We choose some positions. First of all, we take positions $(1,2)$, $\ldots$, $(m-1,m)$. Further, for every $i$, where $1\le i\le m-1$ and $k_i<m$, we choose positions
$$
(i,k_i+1), \ldots, (i,m). \eqno(1)
$$
Now we will do the same for the positions $(i,j)$. And then we fix positions $(2,1),\ldots,(m,m-1)$. For every $j$, where $1\le j\le m-1$ and $k_j<m$, we choose positions
$$
(k_j+1,j), \ldots, (m,j), \eqno(2)
$$
where $k_j$ is the maximal number with $M_{k_j,k_j-1}\cdots M_{j+1,j}\ne 0$. We come to the corresponding facts on elements $c_{ij}$ for $i>j$.

After the work done, we can formulate the following assertion.

\textbf{Proposition 7.1.}\\
\textbf{1.} There exists an isomorphism $\Psi\cong \prod_{(i,j)}C(U(R))$, where pairs $(i,j)$ run over all selected above positions. More precisely, $\Psi\cong\prod_{p}C(U(R))$, where
$$
p=2(m-1)+q,\; q=\sum_{i=1}^{m-1}s_i+\sum_{j=2}^{m}t_j
$$ 
and $s_i$ (resp., $t_j$) is the number of selected positions in $(1)$ (resp., $(2)$).

\textbf{2.} We have isomorphisms
$$
\Psi_0\cong \prod_{m-1}C(U(R)) \text{ and } \Psi/\Psi_0\cong \prod_{(m-1)+q}C(U(R)).
$$
$\lhd$ \textbf{1.} The embedding $\Psi\to \prod_{m^2-m}C(U(R))$ associates an automorphism $\beta\in\Psi$ with the system of invertible central elements
$$
\{c_{ij} \,|\, i,j=1,\ldots,m,\, i\ne j\}
$$
of the ring $R$, where $\beta(y)=c_{ij}y$, $y\in M_{ij}$ (see Section 4 and the above). It follows from the text before the proposition that we can also restrict ourself by elements $c_{ij}$ for pairs $(i,j)$ running over only positions indicated there.

\textbf{2.} Let $\beta\in\Psi_0$ and $c_{ij}$ be the elements from \textbf{1}. For all $i,j,k=1,\ldots,m$, the relation $c_{ik}=c_{ij}\cdot c_{jk}$ holds (see Proposition 4.3 and its proof). From here, we obtain that elements $c_{ij}$ with $i<j$ are products of elements of the form $c_{k,k+1}$. Taking into account that that $1=c_{11}=c_{ij}c_{ji}$, we obtain $c_{ji}= c_{ij}^{-1}$. Therefore, elements $c_{ji}$ with $i<j$ are expressed by elements which are inverse to elements $c_{k,k+1}$. These  considerations lead to the first isomorphism from \textbf{2}. The second isomorphism follows from the first isomorphism and \textbf{1}.~$\rhd$

\textbf{Corollary 7.2, \cite{KryT21}.} If $M^2=0$, then we have isomorphisms
$$
\Psi\cong \prod_{m^2-m}C(U(R)),\; 
\Psi_0\cong \prod_{m-1}C(U(R)),\; \Psi/\Psi_0\cong \prod_{(m-1)^2}C(U(R)).
$$

From Proposition 4.3, Corollary 4.4, and the material of this section, we can formulate when an automorphism from $\Psi$ is inner. We note that automorphisms from $\Psi$ can be called \textsf{multiplicative}.

\textbf{Corollary 7.3.} Let $K$ be a formal matrix algebra over a  ring $R$ with indecomposable factor ring $R/P(R)$. The multiplicative automorphism $\varphi$ is inner if and only if the corresponding to it system of elements
$c_{ij}\in C(U(R))$ ($i,j=1,\ldots,m$) satisfies to the relations $c_{ij}\cdot c_{jk}=c_{ik}$ for all $i,j,k=1,\ldots,m$.

\textbf{Corollary 7.4.} If we add the condition
$$
M_{ij}\cdot M_{jk}\ne 0 \; \text{ for all } i,j,k,
$$
to conditions of Corollary 7.3, then the relation $\Psi=\Psi_0$ holds, i.e., every multiplicative automorphism is inner.

$\lhd$ We have that the $R$-module $M_{ij}\cdot M_{jk}$ is faithful. (This has already been used at the beginning of this section.) This implies relations $c_{ij}\cdot c_{jk}=c_{ik}$.~$\rhd$

\section{Groups $\Omega$ and $\Omega_1$}\label{section8}

As before, $K$ is a formal matrix algebra over the ring $R$. In the beginning of Section 2, the group $\Omega$ was defined; it is the image of the homomorphsm $f\colon \text{Aut}\,K\to \text{Aut}\,L$. In this section, we consider the group $\Omega$ in the case of the ring $K$. The role of this group and the group $\Psi$ has already been mentioned in Section 2 (especially see the end of Section 2). We also define a group $\Omega_1$ as the image of the restriction of the homomorphsm $f$ to $\text{Ker}\,g$.

We return to the decomposition $\text{Aut}\,L=\Gamma\leftthreetimes \Sigma$ from section $6$ and the homomorphsm $g\colon \text{Aut}\,K\to \Sigma$. We recall on the decomposition $K=L\oplus M$. If $M^2=0$, then $\text{Aut}\,K=\text{Ker}\,g\leftthreetimes \Sigma$ \cite [Theorem 13.3]{KryT21}. Next, we have the semidirect decomposition $\Omega=\Omega_1\leftthreetimes \Sigma$ (see the beginning of Section 14 in \cite{KryT21}). If $M^2\ne 0$, then we can only say that $\Omega_1$ is a normal subgroup in $\Omega$ and the factor group $\Omega/\Omega_1$ is isomorpically embedded in the permutation group $\Sigma$. We pay attention to the subgroup $\Omega_1$ especially.

We write several questions on the structure of the group $\Omega_1$.

\textbf{1.} Which automorphisms from $\Gamma$ belong $\Omega_1$?

\textbf{2.} What is the structure of the group $\Omega_1$?

\textsf{Next, we assume that factor ring $R/P(R)$ is indecomposable.}
We will answer the first question and, in one case, the second question. We pay attention to the following circumstance. By Theorem 2.1 and Corollary 5.3, we have the relation
$$
\text{Aut}\,K=\text{In}_1(\text{Aut}\,K)\leftthreetimes \Lambda.
$$
This implies the relation $\text{Ker}\,g=\text{In}_1(\text{Aut}\,K)\leftthreetimes C$, where
$$
C=\text{Ker}\,g\cap\Lambda=
\{\begin{pmatrix}\alpha&0\\ 0&\beta\end{pmatrix}\,|\, \alpha R_i=R_i, \,i=1,\ldots,m\}.
$$
Now we can assert that an automorphism $\alpha\in\Gamma$ is contained in $\Omega_1$ if and only if there is a transformation $\beta$ of the algebra $M$ which is both its automorphism and an isomorphism of $L$-$L$-bimodules $M\to {}_{\alpha}M_{\alpha}$. The last property is equivalent to the property that the matrix $\begin{pmatrix}\alpha&0\\ 0&\beta\end{pmatrix}$ defines automorphism of the algebra $K$ contained in $\text{Ker}\,g$.

Let $\alpha\in\Omega_1$ and let $\varphi=\begin{pmatrix}\alpha&0\\ 0&\beta\end{pmatrix}$ be the corresponding automorphism of the algebra $K$. For any $i,j=1,\ldots,m$, we have the relation
$$
\beta M_{ij}=\beta(e_iMe_j)=\alpha(e_i)\beta(M)\alpha(e_j)=
e_iMe_{j}=M_{ij}.
$$
We set $\beta_{ij}=\beta\big|_{M_{ij}}$ and $\alpha_i=\alpha\big|_{R_i}$. Then $\beta_{ij}\colon M_{ij}\to {}_{\alpha_i}(M_{ij})_{\alpha_j}$ is an isomorphism $R_i$-$R_j$-bimodules (bimodules of the form $_{\alpha}A_{\gamma}$ are defined in Section 1). We can give the following form to the question of which elements of $\Gamma$ are contained in $\Omega_1$. For which automorphisms $\alpha_i\in \text{Aut}\,R_i$ and $\alpha_j\in \text{Aut}\,R_j$ there exist isomorphisms between $R_i$-$R_j$-bimodules $M_{ij}$ and how they are arranged?

We prove a general fact. It generalizes the following result (see \cite[Chapter 2, Proposition 5.2]{Bas68}):

Let $H$ be some algebra and let $\alpha$, $\gamma$ be automorphisms of $H$. There exists an isomorphism of $H$-$H$-bimodules $H\to {}_{\alpha}H_{\gamma}$ if and only if $\alpha^{-1}\gamma$ is an inner automorphism.

Let $k$, $\ell$ be two positive integers and $c=\text{LCM}(k,\ell)$. We denote 
$$
P=M(k,R),\, Q=M(\ell,R),\, H=M(c,R),\, V=M(k\times\ell,R),\, \ell'=\dfrac{c}{k},\, k'=\dfrac{c}{\ell}.
$$
The ring $H$ can be represented as a block matrix ring by two methods: as a block matrix ring over $P$ of order $\ell'$ and as a block matrix ring over $Q$ of order $k'$. It also is a $P$-$Q$-bimodule block matrices over $V$ of size $\ell'\times k'$.

Let $\alpha$ and $\gamma$ be automorphisms of algebras $P$ and $Q$, respectively. They induce the automorphisms of the algebra $H$; we call them \textsf{ring automorphisms}. For them, we leave designations $\alpha$ and $\gamma$, respectively. This agreement is also preserved in the following proposition.

\textbf{Proposition 8.1.} If $\alpha\in\text{Aut}\,P$ and $\gamma\in\text{Aut}\,Q$, then an isomorphism of $P$-$Q$-bimodules $V\to {}_{\alpha}V_{\gamma}$ exists if and only if $\alpha^{-1}\gamma$ is an inner automorphism of the algebra $H$.

$\lhd$ Let we have an isomorphism of $P$-$Q$-bimodules $\beta\colon V\to {}_{\alpha}V_{\gamma}$. The isomorphism $\beta$ induces an $H$-$H$-bimodule isomorphism  
$$
\overline{\beta}\colon H\to {}_{\alpha}H_{\gamma},\;
\overline{\beta}(A)=(\beta(A_{ij}))
$$
for every matrices $A=(A_{ij})\in H$. We mean that the matrix $A$ is represented in the above block form, i.e., $A_{ij}$ are blocks of size $\ell'\times k'$. Consequently, $\alpha^{-1}\gamma$ is an inner automorphism of the algebra $H$.

Now we assume that $\alpha^{-1}\gamma$ is an inner automorphism of the algebra $H$. Consequently, there exists an isomorphism of $H$-$H$-bimodules $\beta\colon H\to {}_{\alpha}H_{\gamma}$.

We take the triangular matrix  algebra $S=\begin{pmatrix}H&H\\ 0&H\end{pmatrix}$. We denote by $\psi$ the automorphism of the algebra $S$ which converts the matrix $\begin{pmatrix}a&c\\0&b\end{pmatrix}$ to the matrix $\begin{pmatrix}\alpha(a)&\beta(c)\\0&\gamma(b)\end{pmatrix}$, i.e., $\psi=\begin{pmatrix}(\alpha,\gamma)&0\\0&\beta\end{pmatrix}$ in the matrix form accepted by us of automorphisms.

Let $e_1,\ldots,e_{\ell'}$ and $f_1,\ldots,f_{k'}$ be diagonal matrix units which correspond to two block partitions of matrices in $H$. We have the relation
$$
\alpha(e_i)=e_i,\; i=1,\ldots,\ell',\; \gamma(f_j)=f_j,\; j=1,\ldots,k'.
$$
This implies that $\psi$ induces the automorphism of triangular matrix algebras $\begin{pmatrix}e_iHe_i&e_iHf_j\\ 0&f_jHf_j\end{pmatrix}$; in fact, this means that $\psi$ induces the automorphism of the algebra $\begin{pmatrix}P&V\\ 0&Q\end{pmatrix}$. Consequently, $\beta\big|_V$ is isomorphism $P$-$Q$-bimodules $V\to {}_{\alpha}V_{\gamma}$.~$\rhd$

Let $n_i$ be the order of matrices in the ring $R_i$, $i=1,\ldots,m$. We set $c_{ij}=\text{LCM}(n_i,n_j)$ for all pairwise distinct $i,j=1,\ldots,m$. We denote by $H_{ij}$ the matrix ring $M(c_{ij},R)$. It is a ring block matrices over $R_i$ and over $R_j$ and also is an $R_i$-$R_j$-bimodule of block matrices over $M_{ij}$. We assume that automorphisms of the rings $R_i$ and $R_j$ are rings automorphisms of the algebra $H_{ij}$.

\textbf{Theorem 8.2.} The automorphism $\alpha=(\alpha_i)$ of the algebra $L=\oplus_{i=1}^mR_i$ belongs to the group $\Omega_1$ if and only if $\alpha_i^{-1}\alpha_j$ is an inner automorphism of the algebra $H_{ij}$ for all distinct subscripts $i$ and $j$. 

$\lhd$ \textsf{Necessity.} Let $\alpha\in\Omega_1$. Consequently, there exists an automorphism $\varphi=\begin{pmatrix}\alpha&0\\ 0&\beta\end{pmatrix}$ of the algebra $K$, where $\beta$ is an automorphism $L$-$L$-bimodules $M\to {}_{\alpha}M_{\alpha}$. The restricition $\beta$ to $M_{ij}$ is an $R_i$-$R_j$-bimodule isomorphism $M_{ij}\to {}_{\alpha_i}(M_{ij})_{\alpha_j}$. According to Proposition 8.1, $\alpha_i^{-1}\alpha_j\in \text{In(Aut}\,H_{ij})$.

\textsf{Sufficiency.} By Proposition 8.1, there is an isomorphism
of $R_i$-$R_j$-bimodules $\beta_{ij}\colon M_{ij}\to {}_{\alpha_i}(M_{ij})_{\alpha_j}$, $i,j=1,\ldots,m$, $i\ne j$. Let $\beta=\sum_{i,j=1,i\ne j}^m\beta_{ij}$. Then $\beta$ is an isomorphism of $L$-$L$-bimodules $M\to {}_{\alpha}M_{\alpha}$. For the transformation $\begin{pmatrix}\alpha&0\\ 0&\beta\end{pmatrix}$ of the algebra $K$ to be its automorphism, it suffices to verify that the equality
$$
\beta_{ik}(ab)=\beta_{ij}(a)\beta_{jk}(b) \eqno (*)
$$
holds for any elements $a\in M_{ij}$, $b\in M_{jk}$ and of all pairwise distinct subscripts $i,j,k$.

We fix three mentioned subscripts $i,j,k$ and define another matrix ring. We set $H=M(d,R)$, where $d=\text{LCM}(n_i,n_j,n_k)$. This ring $H$ is a block matrix ring over the bimodules $M_{ij}$, $M_{jk}$, $M_{ik}$.

We assume that automorphisms $\alpha_i$, $\alpha_j$, $\alpha_k$ are ring automorphisms of the algebra $H$. We consider isomorphisms $\beta_{ij}$, $\beta_{jk}$, $\beta_{ik}$ as bimodule isomorphisms
$$
H\to {}_{\alpha_i}H_{\alpha_j},\; H\to {}_{\alpha_j}H_{\alpha_k},\; H\to {}_{\alpha_i}H_{\alpha_k}
$$
respectively. At the same time, products $\alpha_i^{-1}\alpha_j$, $\alpha_j^{-1}\alpha_k$, $\alpha_i^{-1}\alpha_k$ are inner automorphisms of the algebra $H$. They induce the same bimodule isomorphisms $\beta_{ij}$, $\beta_{jk}$ and $\beta_{ik}$ as above.

These bimodule isomorphisms act as follows (see \cite[Section2]{KryT21}). There exist invertible elements $u,v,w\in H$ such that we have the relation
$$
\beta_{ij}(x)=\alpha_i(x)\alpha_i(u)=\alpha_i(u)\alpha_j(x),
$$
$$
\beta_{jk}(y)=\alpha_j(y)\alpha_j(v)=\alpha_j(v)\alpha_k(y),
$$
$$
\beta_{ik}(z)=\alpha_i(z)\alpha_i(w)=\alpha_i(w)\alpha_k(z)
$$
for all $x,y,z\in H$.

We verify that the relation
$$
\beta_{ik}(xy)=\beta_{ij}(x)\beta_{jk}(y) \eqno (**)
$$
holds in $H$ for any $x,y\in H$. It follows from the relation
$$
\alpha_i^{-1}\alpha_k=\alpha_i^{-1}\alpha_j\alpha_j^{-1}\alpha_k
$$
that the relation $w=vu$ holds (it is also necessary to take into account that the elements $u,v,w$  are defined up to invertible central elements). Now we have the relation
$$
\beta_{ik}(xy)=\alpha_i(xy)\alpha_i(w)=
\alpha_i(x)\alpha_i(y)\alpha_i(v)\alpha_i(u).
$$
Also we have the relation
$$
\beta_{ij}(x)\beta_{jk}(y)=\alpha_i(x)\alpha_i(u)\alpha_j(y)\alpha_j(v)=
$$
$$
=\alpha_i(x)\alpha_i(y)\alpha_i(u)\alpha_j(v)=
\alpha_i(x)\alpha_i(y)\alpha_i(v)\alpha_i(u).
$$
Thus, the relation $(**)$ is proved. In it, matrix multiplication is executed in block form. This implies the relation $(*)$.~$\rhd$

In the remaining part of the section, we touch on the structure problem for the group $\Omega_1$. This problem seems to be quite complicated. We will find the structure of the group $\Omega_1$ under one condition on the integers $n_i$. We recall that $n_i$ is the order of matrices in the ring $R_i$.

Let $n_s$ be the least of integers $n_1,\ldots,n_m$. Next, we assume that the algebra $K$ satisfies the following property: $n_s$ divides each of the integers $n_i$. Under such an assumption, every ring $R_i$ is a block matrix ring over $R_s$ of order $n_i/n_s$. Therefore, every automorphism $\alpha_s$ of the algebra $R_s$ is extended to ring automorphism $\alpha_i$ of the algebra $R_i$. We call the obtained automorphism $(\alpha_1,\alpha_2,\ldots,\alpha_m)$ of $L$ a \textsf{scalar automorphism}. Also we write it in the form $(\alpha_s,\alpha_s,\ldots,\alpha_s)$, i.e., we identify $\alpha_i$ with $\alpha_s$.

We denote by $D$ the subgroup of all scalar automorphisms of the algebra $L$. The following result extends \cite[Corollary 14.3]{KryT21}.

\textbf{Corollary 8.3.}\\
\textbf{1.} The relation $\Omega_1=\text{In(Aut}\,L)\cdot D$ holds.

\textbf{2.} There exists an isomorphism $\Omega_1/\text{In(Aut}\,L)\cong \text{Out}\,R_s$.

$\lhd$ \textbf{1.} Let $\alpha=(\alpha_i)\in \Omega_1$. For every $j=1,\ldots,m$, we denote by $\mu_j$ the automorphism $\alpha_j\alpha_s^{-1}$ of the algebra $R_j$. It follows from Theorem 8.2 that $\mu_j$ is an inner automorphism of the algebra $R_j$. Now it follows from the relations $\alpha_j=\mu_j\alpha_s$ ($j=1,\ldots,m$) that
$$
\alpha=(\mu_1,\ldots,\mu_m)(\alpha_s,\ldots,\alpha_s)=\mu\gamma\in\text{In(Aut}\,L)\cdot D.
$$
The inclusion $\Omega_1\subseteq \text{In(Aut}\,L)\cdot D$ is proved.

We prove the converse inclusion. The inclusion $\text{In(Aut}\,L)\subseteq \Omega_1$ holds always. Now we take an arbitrary automorphism $\gamma=(\alpha,\ldots,\alpha)$ from $D$, where $\alpha\in R_s$. The inclusion $\gamma\in\Omega_1$ follows from Theorem 8.2.

\textbf{2.} The assertion follows from \textbf{1}.~$\rhd$

\section[Automorphisms of Triangular Matrix Rings]{Automorphisms of Triangular Matrix Rings}\label{section9}

In this section, $R$ is an arbitrary algebra over some commutative ring $T$. The (upper) triangular matrix ring over $R$ is denoted by $T(n,R)$. We denote this matrix ring by $K$. We consider the ring $K$ as a splitting extension: $T(n,R)=K=L\oplus M$, where the symbols $L$ and $M$ have the same meaning. It also is convenient to assume that the ring $L$ is the sum $R_1\oplus \ldots\oplus R_n$ and the $L$-$L$-bimodule $M$ is equal to $\oplus_{i,j=1, i<j}^nM_{ij}$.

Every automorphism $\alpha$ of the algebra $R$ induces the automorphism $\overline{\alpha}$ of the algebra $K$, where
$$
\overline{\alpha}(a_{ij})=(\alpha(a_{ij})),\;(a_{ij})\in K.
$$
The automorphism $\overline{\alpha}$ is called the \textsf{ring automorphism induced by $\alpha$.} Such automorphisms were considered in the previous section.

The theorem below is another formulation of one theorem from \cite{Kop96}; the proof differs from the proof in \cite{Kop96}.

\textbf{Theorem 9.1.} Every triangular automorphism of the algebra $T(n,R)$ is a product of an inner automorphism and ring automorphism.

$\lhd$ Let $\varphi$ be some triangular automorphism of the algebra $K$. We represent it as a product of an inner and a diagonal automorphisms, since Theorem 2.1(1) is true for the subgroup of triangular automorphisms. Therefore, we can assume that $\varphi$ is a diagonal automorphism. Thus, $\varphi=\begin{pmatrix}\alpha&0\\ 0&\beta\end{pmatrix}$, where $\alpha$ is an automorphism of the algebra $L$, $\beta$ is an automorphism of the algebra $M$ and  an isomorphism of $L$-$L$-bimodules $M\to {}_{\alpha}M_{\alpha}$.

Since $e_1,\ldots,e_n$ is a complete orthogonal system of central idempotents of the ring $L$, we have that the system $\alpha(e_1),\ldots,\alpha(e_n)$ satisfies the same properties. For every $i$, we write
$$
\alpha(e_i)=f_1^{(i)}+f_2^{(i)}+\ldots+f_n^{(i)},\eqno (1)
$$
where $f_k^{(i)}\in R_k$, $k=1,\ldots,n$. Elements $f_k^{(1)},f_k^{(2)},\ldots, f_k^{(n)}$ form a complete orthogonal system of central idempotents of the ring $R_k$.

There exist inclusions
$$
M\supset M^2\supset\ldots\supset M^{n-1}\supset M^n=0,
$$
where $M^{n-1}=M_{1n}$. Therefore, the relation $\varphi M_{1n}=M_{1n}$ is true. The relation $e_1Ke_n=M_{1n}$ implies the relation
$$
\varphi(e_1Ke_n)=\alpha(e_1)K\alpha(e_n)=\varphi M_{1n}=M_{1n}=R.
$$
Next, we can write the relation
$$
\alpha(e_1)K\alpha(e_n)=\alpha(e_1)(L\oplus M)\alpha(e_n)=f_1^{(1)}Mf_n^{(n)}.
$$
Thus, $f_1^{(1)}Mf_n^{(n)}=R$. Consequently, $f_1^{(1)}=1$ (i.e., $e_1$), $f_n^{(n)}=1$ (i.e., $e_n$) and
$$
f_1^{(2)}=\ldots =f_1^{(n)}=0,\;
f_n^{(1)}=\ldots =f_n^{(n-1)}=0.
$$
We use the induction on $n$ to show that $\alpha(e_i)=e_i$, $i=1,\ldots,n$. It's already been proven for $n=2$. Let $n\ge 3$. We represent the ring $K$ in the form of a block-triangular matrix ring of order 2:
$$
K=\begin{pmatrix}S&N\\ 0&R_n\end{pmatrix}, \text{ where}
$$
$$
S=\begin{pmatrix}
R_1&M_{12}&\ldots&M_{1,n-1}\\
0&R_2&\ldots&M_{2,n-1}\\
\ldots&\ldots&\ldots&\ldots\\
0&0&\ldots&R_{n-1}
\end{pmatrix},
N=\begin{pmatrix}
M_{1n}\\
M_{2n}\\
\ldots\\
M_{n-1,n}
\end{pmatrix}.
$$
We verify that $\varphi S=S$, i.e., $\varphi$ induces the automorphism of the ring $S$. First, it follows from relations $(1)$ that
$$
\varphi(R_1),\ldots,\varphi(R_{n-1})\subseteq
R_1\oplus\ldots\oplus R_{n-1}.
$$
Let $M_{ik}$ be an arbitrary bimodule, where $i<k$ and $k\ne n$. Since the element $\alpha(e_k)$ has the zero component in $R_n$, we have that the last column of all matrices from $\alpha(e_i)M\alpha(e_k)$ consists of zeros. Therefore, $\varphi(M_{ik})\subseteq S$ and $\varphi(S)\subseteq S$. Similarly, we obtain $\varphi^{-1}(S)\subseteq S$.

Let $\begin{pmatrix}\alpha_1&\gamma_1\\ 0&\beta_1\end{pmatrix}$ and $\begin{pmatrix}\alpha_2&\gamma_2\\ 0&\beta_2\end{pmatrix}$ be matrices of automorphisms $\varphi$ and $\varphi^{-1}$, respectively with respect to the decomposition $K=L_1\oplus N$, where $L_1=S\oplus R_n$. Then $\alpha_1\alpha_2=1=\alpha_2\alpha_1$, i.e., $\alpha_1=\varphi\big|_S$ is an automorphism of the algebra $S$. By the induction hypothesis, we have $f_2^{(2)}=\ldots =f_{n-1}^{(n-1)}=1$ and all <<remaining>> $f_i^{(j)}$ are equal to zero. It is proved that
$$
\alpha(e_1)=e_1,\ldots,\alpha(e_n)=e_n.
$$
We return to the initial representation of the ring $K$; namely, to the relation 
$K=L\oplus M$. We take an arbitrary subscript $i$ and the element $x\in R_i$. Then
$$
x=xe_i,\; \varphi(x)=\varphi(x)\varphi(e_i)=\varphi(x)e_i\in L\cap Ke_i=R_i.
$$
Therefore, $\alpha R_i=R_i$ for all $i=1,\ldots,n$.

We denote by $\alpha_i$ the restriction of $\alpha$ to $R_i$ ($i=1,\ldots,n$). Then $\alpha_i$ is an automorphism of the ring $R_i$. We also have the relation
$$
\beta M_{ij}=\beta(e_iMe_j)=\alpha(e_i)M\alpha(e_j)=e_iMe_j=M_{ij}. 
$$
We denote by $\beta_{ij}$ the restriction $\beta\big|_{M_{ij}}$. Since $\beta_{ij}(xyz)=\alpha_i(x)\beta(y)\alpha_j(z)$ for any $x\in R_i$, $z\in R_j$, $y\in M_{ij}$, we have that $\beta_{ij}$ is an isomorphism of $L$-$L$-bimodules $M_{ij}$ and ${}_{\alpha_i}(M_{ij})_{\alpha_j}$. In such a case, $\alpha_i\alpha_j^{-1}$ is an inner automorphism of the algebra $R$ \cite[Part 2, Proposition 5.2]{Bas68}. For every $i=1,\ldots,n$, we set $\mu_i=\alpha_i\alpha_1^{-1}$. Then
$$
\alpha_i=\mu_i\alpha_1, \; \alpha=(\mu_1,\ldots,\mu_n)(\alpha_1,\ldots,\alpha_1),
$$
where $\mu_1,\ldots,\mu_n$ are inner automorphisms of the ring $R$. More briefly, $\alpha=\mu\gamma$, where $\mu$ is an inner automorphism and $\gamma$ is a ring automorphism of the ring $L$.

Automorphisms $\mu$ and $\gamma$ induce an inner automorphism and a ring automorphism of the ring $K$. Let us keep the notation $\mu$ and $\gamma$, respectively. We set $\zeta=(\mu\gamma)^{-1}\varphi$ and show that $\zeta$ is an inner automorphism. With respect to the decomposition $K=L\oplus M$,  its matrix is of the form $\varphi=\begin{pmatrix}1&0\\ 0&\rho\end{pmatrix}$, where $\rho$ is an automorphism of the algebra $M$ and an automorphism of the $L$-$L$-bimodule $M$. We set $\rho_{ij}=\rho\big|_{M_{ij}}$ for all $i,j$ with $i<j$. Then $\rho_{ij}$ is an automorphism of the $R$-$R$-bimodule $M_{ij}$, i.e., $\rho_{ij}$ is an automorphism of the $R$-$R$-bimodule $R$, since 
$$
\rho M_{ij}=\rho(e_iMe_j)=e_iMe_j=M_{ij}.
$$
Consequently, $\rho_{ij}$ acts on $M_{ij}$ as a multiplication by some invertible central element $c_{ij}$ of the ring $R$.

Thus, the automorphism $\zeta$ defines a system of invertible central elements $c_{ij}$ ($i,j=1,\ldots,n$, $i<j$). At the same time, the relation $c_{ij}\cdot c_{jk}=c_{ik}$ hold for all $i,j,k=1,\ldots,n$ such that $i<j<k$. To verify that $\zeta$ is an inner automorphism, we can use the argument which is similar to the argument from the proof of Proposition 4.3.

As a result, we obtain the relation $\varphi=\mu\gamma\zeta$ in which $\mu,\zeta$ are inner automorphisms and $\gamma$ is a ring automorphism. By reformulating, we obtain the relation $\varphi=\mu\zeta'\gamma$, where $\zeta'=\gamma\zeta\gamma^{-1}$ is an inner automorphism. Thus, the automorphism $\varphi$ is equal to a product of an inner automorphism and a ring automorphism, which is required.~$\rhd$

The following  Corollaries 9.2 and 9.3, directly follow from the above theorem, Proposition 3.2 and Lemma 3.3.

\textbf{Corollary 9.2.} Any automorphism of the algebra $T(n,R)$ is a product of an inner automorphism and a ring automorphism in each of the following cases.

\textbf{1)} $R$ is a strongly indecomposable algebra.

\textbf{2)} $R$ is a semiprime or normal algebra, \cite{Jon95}.

\textbf{Corollary 9.3, \cite{Kez90}.} If $R$ is a commutative ring, then all automorphisms of the $R$-algebra $T(n,R)$ are inner.

\section{Group $\text{Aut}\,K$ with $K=M(n,R)$}\label{section10}

Similar to the previous section, $R$ is an algebra over a commutative ring $T$. We give some remarks on automorphisms of the $T$-algebra $K=M(n,R)$.

In previous sections, we used methods based on splitting extensions. However, it is not applicable to the algebra $M(n,R)$.
One of the possible approaches to the study the algebra $\text{Aut}\,K$ is based on interrelations of this group with the Picard group of the ring $R$. Such approach is used in \cite{RosZ61} for separable algebras and in \cite{Isa80} for the $T$-algebra $M(n,T)$.

The Picard group $\text{Pic}\,S$ of the $T$-algebra $S$ is defined as the class group $[P]$ of isomorphic invertible $S$-$S$-bimodules $P$ with operation $[P]\cdot [Q]=[P\otimes_SQ]$; see \cite{Bas68}.

As usual, finitely generated projective generators are called \textsf{progenerators}.
The direct sum of $n$ copies of the module $A$ is denoted by $A^n$. We use the following fact \cite[Part 2, Proposition 5.2]{Bas68}.

Let $A$ be an $R$-$R$-bimodule and let we have an isomorphism $A\cong R$ of left $R$-modules. Then for some $\alpha\in \text{Aut}\,R$, there exists an isomorphism of $R$-$R$-bimodules $A\cong {}_1R_{\alpha}$.

We give a familiar result \cite[Part 2, Proposition 5.3]{Bas68}. 

\textbf{Proposition 10.1.} Let $Q$ be a left progenerator 
$R$-module and let the endomorphism ring $S=\text{End}_RQ$ be considered as a $T$-algebra. Then we have a group exact sequence 
$$
1\to\text{In(Aut}_TS)\to\text{Aut}_TS\stackrel{\eta_Q}{\longrightarrow} \text{Pic}\,R,
$$
$$
\text{where Im}\,\eta_Q=\{[P]\in\text{Pic}\,R\,|\,P\otimes_RQ\cong Q \text{ as left } R\text{-modules}\}.
$$
We apply Proposition 10.1 to the algebra $K=M(n,R)$. We take $R^n$ as the module $Q$. Then we can identify the algebra $S$ with $K$. Next, if $[P]\in\text{Pic}\,R$, then there exists an isomorphism of left 
$R$-modules $P\otimes_RQ\cong Q$ if and only if left $R$-modules $P^n$ and $R^n$ are isomorphic. We can write the following result.

\textbf{Proposition 10.2.} There exists an exact sequence of groups
$$
1\to\text{In(Aut}\,K)\to\text{Aut}\,K\stackrel{\eta}{\longrightarrow} \text{Pic}\,R,\eqno (1)
$$
$$
\text{where Im}\,\eta=\{[P]\in\text{Pic}\,R\,|\,P^n\cong R^n \text{ as left } R\text{-modules}\}.
$$
The homomorphism $\eta$ is the composition of homomorphsms
$$
\text{Aut}\,K\stackrel{\zeta}{\longrightarrow} \text{Pic}\,K\stackrel{\theta}{\longrightarrow} \text{Pic}\,R.
$$
Here $\zeta\colon \varphi\to[_1K_{\varphi}]$; in addition, $\theta$ is the canonical isomorphism which exists, since categories $K$-mod and $R$-mod are equivalent. More precisely, we identify $R^n$ with the module of row vectors of length $n$ and denote the $R$-$K$-bimodule $R^n$ by $M$. We also denote the $K$-$R$-bimodule of column vectors $R^n$ by $N$. Then
$$
\theta\colon [Q]\to [M\otimes_KQ\otimes_KN],\;
\theta^{-1}\colon [P]\to [N\otimes_RP\otimes_RM].
$$
Thus, we have
$$
\eta\colon \varphi\to [M\otimes_K {}_1K_{\varphi}\otimes_KN]=
[M_{\varphi}\otimes_KN].
$$

\textbf{Proposition 10.3.} The following two assertions are equivalent.

\textbf{1)} Every automorphism of the algebra $K=M(n,R)$ is a product of an inner automorphism and a ring automorphism.

\textbf{2)} If $[P]\in \text{Pic}\,R$ and left $R$-modules $P^n$ and $R^n$ are isomorphic, then left $R$-modules $P$ and $R$ are isomorphic.

$\lhd$ \textbf{1)\,$\Rightarrow$\,2).} Let $[P]\in\text{Pic}\,R$ and left $R$-modules $P^n$ and $R^n$ are isomorphic. According to Proposition 10.2, we have $[P]\in\text{Im}\,\eta$. We choose $\varphi\in\text{Aut}\,K$ such that $\eta(\varphi)=[P]$. Based on \textbf{1)}, we have $\varphi=\mu\psi$, where $\mu$ is an inner automorphism and $\psi$ is a ring automorphism. Therefore, $\eta(\varphi)=\eta(\psi)=[P]$. Let $\psi$ be defined by an automorphism $\alpha\in \text{Aut}\,R$. Then we have isomorphisms
$$
P\cong M_{\psi}\otimes_KN\cong {}_1R_{\alpha}.
$$
\textbf{2)\,$\Rightarrow$\,1).} If $\varphi$ is an arbitrary automorphism of the algebra $K$, then 
$$
\eta(\varphi)=[P]\in\text{Pic}\,R,\; P^n\cong R^n \text{ as left } R\text{-modules}.
$$
According to \textbf{2)}, there exists an isomorphism ${}_RP\cong {}_RR$.
Consequently, ${}_1P_1\cong {}_1P_{\alpha}$ for some $\alpha\in\text{Aut}\,R$. We also have the following relation in the group $\text{Pic}\,K$:
$$
[N\otimes_R {}_1R_{\alpha}\otimes_RM]=[N_{\alpha}\otimes_RM]=[_1K_{\overline{\alpha}}].
$$
Now we can write the following relation:
$$
\varphi=\eta^{-1}([P])=\eta^{-1}([_1R_{\alpha}])=
\zeta^{-1}\theta^{-1}([_1R_{\alpha}])=
$$
$$
=\zeta^{-1}([N\otimes_R {}_1R_{\alpha}\otimes_RM])=
\zeta^{-1}([_1K_{\overline{\alpha}}]).
$$
In other words, $\zeta(\varphi)=[_1K_{\overline{\alpha}}]=\zeta(\overline{\alpha})$,
where $\overline{\alpha}$ is the automorphism induced by $\alpha$.
Consequently, $\varphi=\mu\overline{\alpha}$, where $\mu$ is some inner automorphism of the algebra $K$.~$\rhd$

\textbf{Corollary 10.4.} Let the ring $R$ do not have non-trivial idempotents and every left $R$-progenerator satisfies the isomorphism property of direct decompositions. Then any automorphism of the algebra $M(n,R)$ is equal to a product of an inner automorphism and ring automorphism.

$\lhd$ Let $[P]\in \text{Pic}\,R$ and let left $R$-modules $P^n$, $R^n$ be isomorphic. It follows from the isomorphism $\text{End}_RP\cong R$ that $_RP$ is an indecomposable module. Therefore, $_RP\cong {}_RR$ and we can use Proposition 10.3.~$\rhd$

\textbf{Corollary 10.5.} If $R$ is a local ring or a principal left ideal domain, then any automorphism of the algebra $M(n,R)$ is equal to a product of an inner automorphism and a ring automorphism.

$\lhd$ Let left $R$-modules $P^n$ and $R^n$ be isomorphic. Then $\text{End}_RP\cong R$. If the ring $R$ is local, then we immediately obtain $_RP\cong {}_RR$. If $R$ is a principal left ideal domain, then the module $_RP$ is isomorphic to a direct sum of left ideals of the ring $R$. Therefore, $_RP\cong {}_RR$ and we can use Proposition 10.3.~$\rhd$

\textbf{Remark 10.6.} Under conditions of Corollaries 10.4 and 10.5, we have an isomorphism $\text{Out}\,K\cong\text{Out}\,R$.

The Tuganbaev's study is supported by Russian Scientific Foundation. 

\addtocontents{toc}{\textbf{$\quad\;$Bibliography$\qquad\qquad\qquad\qquad\qquad\qquad\qquad\qquad\qquad\qquad$ \pageref{biblio}}\par}

\label{biblio}

\end{document}